\documentclass[12pt]{amsart}
\usepackage{mathrsfs}
\usepackage{amssymb}
\usepackage{verbatim}
\usepackage{hyperref}
\usepackage{color}
\usepackage{amsfonts}
\usepackage{mathrsfs}
\usepackage{amsmath}
\usepackage{amssymb}
\usepackage{bbm}

\begin{document}
\newtheorem{theorem}{Theorem}
\newtheorem{proposition}[theorem]{Proposition}
\newtheorem{conjecture}[theorem]{Conjecture}
\def\theconjecture{\unskip}
\newtheorem{corollary}[theorem]{Corollary}
\newtheorem{lemma}[theorem]{Lemma}
\newtheorem{sublemma}[theorem]{Sublemma}
\newtheorem{observation}[theorem]{Observation}
\theoremstyle{definition}
\newtheorem{definition}{Definition}
\newtheorem{notation}[definition]{Notation}
\newtheorem{remark}[definition]{Remark}
\newtheorem{question}[definition]{Question}
\newtheorem{questions}[definition]{Questions}
\newtheorem{example}[definition]{Example}
\newtheorem{problem}[definition]{Problem}
\newtheorem{exercise}[definition]{Exercise}

\numberwithin{theorem}{section}
\numberwithin{definition}{section}
\numberwithin{equation}{section}

\def\earrow{{\mathbf e}}
\def\rarrow{{\mathbf r}}
\def\uarrow{{\mathbf u}}
\def\varrow{{\mathbf V}}
\def\tpar{T_{\rm par}}
\def\apar{A_{\rm par}}

\def\reals{{\mathbb R}}
\def\torus{{\mathbb T}}
\def\heis{{\mathbb H}}
\def\integers{{\mathbb Z}}
\def\naturals{{\mathbb N}}
\def\complex{{\mathbb C}\/}
\def\distance{\operatorname{distance}\,}
\def\support{\operatorname{support}\,}
\def\dist{\operatorname{dist}\,}
\def\Span{\operatorname{span}\,}
\def\degree{\operatorname{degree}\,}
\def\kernel{\operatorname{kernel}\,}
\def\dim{\operatorname{dim}\,}
\def\codim{\operatorname{codim}}
\def\trace{\operatorname{trace\,}}
\def\Span{\operatorname{span}\,}
\def\dimension{\operatorname{dimension}\,}
\def\codimension{\operatorname{codimension}\,}
\def\nullspace{\scriptk}
\def\kernel{\operatorname{Ker}}
\def\ZZ{ {\mathbb Z} }
\def\p{\partial}
\def\rp{{ ^{-1} }}
\def\Re{\operatorname{Re\,} }
\def\Im{\operatorname{Im\,} }
\def\ov{\overline}
\def\eps{\varepsilon}
\def\lt{L^2}
\def\diver{\operatorname{div}}
\def\curl{\operatorname{curl}}
\def\etta{\eta}
\newcommand{\norm}[1]{ \|  #1 \|}
\def\expect{\mathbb E}
\def\bull{$\bullet$\ }

\def\xone{x_1}
\def\xtwo{x_2}
\def\xq{x_2+x_1^2}
\newcommand{\abr}[1]{ \langle  #1 \rangle}
\newcommand{\Norm}[1]{ \left\|  #1 \right\| }
\newcommand{\set}[1]{ \left\{ #1 \right\} }
\def\one{\mathbf 1}
\def\whole{\mathbf V}
\newcommand{\modulo}[2]{[#1]_{#2}}
\def \essinf{\mathop{\rm essinf}}
\def\scriptf{{\mathcal F}}
\def\scriptg{{\mathcal G}}
\def\scriptm{{\mathcal M}}
\def\scriptb{{\mathcal B}}
\def\scriptc{{\mathcal C}}
\def\scriptt{{\mathcal T}}
\def\scripti{{\mathcal I}}
\def\scripte{{\mathcal E}}
\def\scriptv{{\mathcal V}}
\def\scriptw{{\mathcal W}}
\def\scriptu{{\mathcal U}}
\def\scriptS{{\mathcal S}}
\def\scripta{{\mathcal A}}
\def\scriptr{{\mathcal R}}
\def\scripto{{\mathcal O}}
\def\scripth{{\mathcal H}}
\def\scriptd{{\mathcal D}}
\def\scriptl{{\mathcal L}}
\def\scriptn{{\mathcal N}}
\def\scriptp{{\mathcal P}}
\def\scriptk{{\mathcal K}}
\def\frakv{{\mathfrak V}}
\def\C{\mathbb{C}}
\def\D{\mathcal{D}}
\def\I{\mathcal{I}}
\def\R{\mathbb{R}}
\def\Rn{{\mathbb{R}^n}}
\def\Sn{{{S}^{n-1}}}
\def\M{\mathcal{M}}
\def\N{\mathbb{N}}
\def\Q{{\mathcal{Q}}}
\def\Z{\mathbb{Z}}
\def\F{\mathcal{F}}
\def\L{\mathcal{L}}
\def\S{\mathcal{S}}
\def\T{\mathcal{T}}
\def\supp{\operatorname{supp}}
\def\dist{\operatorname{dist}}
\def\essi{\operatornamewithlimits{ess\,inf}}
\def\esss{\operatornamewithlimits{ess\,sup}}

\author{Mingming Cao}\address{
School of Mathematical Sciences\\
Beijing Normal University\\
Laboratory of Mathematics and Complex Systems\\
Ministry of Education\\
Beijing 100875\\
People's Republic of China}\email{m.cao@mail.bnu.edu.cn}
\author{Qingying Xue}

\address{
School of Mathematical Sciences\\
Beijing Normal University\\
Laboratory of Mathematics and Complex Systems\\
Ministry of Education\\
Beijing 100875\\
People's Republic of China}

\email{qyxue@bnu.edu.cn}

\thanks{The second author was supported partly by NSFC
(No. 11471041), the Fundamental Research Funds for the Central Universities (No. 2014KJJCA10) and NCET-13-0065.\\
\indent Corresponding author: Qingying Xue \indent Email: qyxue@bnu.edu.cn}


\date{February 20, 2016.}
\title[Two-weight Boundedness]{Two-weight entropy Boundedness of Multilinear Fractional Type Operators}
\maketitle

\begin{abstract}
This paper will be devoted to study the two-weight norm inequalities of the multilinear fractional maximal operator $\M_{\alpha}$ and the multilinear fractional integral operator $\I_{\alpha}$. The entropy conditions in the multilinear setting will be introduced and the entropy bounds for $\M_\alpha$ and $\I_\alpha$ will be given. 
\end{abstract}

\section{Introduction}
\subsection{Background}
Let $M_{\alpha}$ and $I_\alpha$ be the fractional maximal operator and fractional integral operator defined by
$$
M_\alpha f(x):= \sup_{Q} |Q|^{\frac{\alpha}{n}} \langle f \rangle_Q \cdot \mathbf{1}_Q(x), \ \ I_{\alpha}f(x) := \int_{\Rn} \frac{f(y)}{|x-y|^{n-\alpha}} dy,\ \ \ 0 \leq \alpha < n.
$$
In 1982, Saywer first \cite{Saywer1} showed that
$M_{\alpha}(\cdot \sigma): L^p(\sigma) \rightarrow L^q(w)$ holds if and only if $(w,\sigma)$ satisfies the following testing condition
$$
[w,\sigma]_{S_{(p,q)}} := \sup_{Q} \sigma(Q)^{-1/p} \big\| \mathbf{1}_Q M_{\alpha}(\mathbf{1}_Q \sigma) \big\|_{L^q(w)} < \infty.
$$
Subsequently, using the similar testing conditions, Saywer \cite{Saywer2, Saywer3} gave some characterizations of two weighted weak and strong type inequalities of $I_\alpha$.

After the works of Saywer, many works have been done in the characterizations of two weighed boundedness of continuous operators. Among such achievements are the celebrated works of
Hyt\"{o}nen \cite{H},  Lacey \cite{Lacey1, Lacey2}, Lacey et al \cite{LSTS}, which demonstrated the characterizations of the two weighted $L^2$
inequality of Hilbert transform in terms of Saywer type testing conditions and two weighted $A_2$ condition. Recently, Lacey and Li \cite{LL-g} gave a characterization of two-weight norm inequalities for the classical Littlewood-Paley $g$-function. Still more recently, Cao, Li and Xue \cite{CLX} obtained the characterization of two weighted inequalities for the $g_{\lambda}^*$-function with more general fractional type of Poisson kernels. As for the discrete operators, on the one hand, two-weight characterizations of martingale transforms and dyadic shifts were presented by Nazarov et al \cite{NTV} and Hyt\"{o}nen \cite{H-A2}. The two weighted $L^p(\sigma) \rightarrow L^q(w)$-type inequalities of positive dyadic operators were established by Nazarov et al \cite{NTV} with $p=q=2$, Lacey et al \cite{LST} with $p < q$ and Hyt\"{o}nen \cite{HHL} with $p,q \in(1,\infty)$. On the other hand, in order to study the sufficient condition for the two weight inequalities of the singular integral operators, Treil and Volberg \cite{TV} introduced the entropy conditions. Later on, the entropy conditions were used to obtain the two weight norm inequalities of intrinsic square functions and fractional maximal and integral operators by Lacey, Li \cite{LL-g} and Rahm, Spencer \cite{RS}, respectively.

In the multilinear setting, several works also have already been done for the mulitilinear fractional maximal operator $\M_\alpha$ and fractional integral operators $\I_\alpha$ ($0\leq \alpha< mn$), which are defined by
\begin{eqnarray*}
\M_\alpha(\vec{f})(x)
=\sup_{Q}|Q|^{\frac{\alpha}{n}} \prod_{i=1}^m \langle |f_i| \rangle_Q  \cdot \mathbf{1}_Q(x), \ \
\I_\alpha(\vec{f})(x)= \int_{(\Rn)^m} \frac{\prod_{i=1}^m f_i(x-y_i)}{|(y_1,\ldots,y_m)|^{mn-\alpha}}d\vec{y}.
\end{eqnarray*}
In 2013, Chen and Dami\'{a}n \cite{CD} first gave some sufficient conditions for the two-weight inequalities of the multilinear maximal operator $\mathcal{M}_0$. In 2015, Li and Sun \cite{LS-1} considered the problem of two weighted inequalities of multilinear fractional maximal operator $\mathcal{M}_{\alpha}$. However, their method is not valid for the case $0 \leq \alpha < n(1/p-1/{ \max\{ p_i \} })$. In 2016, Cao and Xue \cite{CX} extended the ranges of exponents to $0 \leq \alpha < mn$ by applying the atomic decomposition of tent space. Moreover, Cao, Xue and Yabuta \cite{CXY} defined and studied the multilinear fractional strong maximal operator and the corresponding multiple weights associated with rectangles. Under the dyadic reverse doubling condition, a necessary and sufficient condition for two-weight inequalities of the multilinear fractional strong maximal operator was given.

It is well known that it is difficult to give a two-weight characterization of $\M_\alpha$ and $\I_\alpha$ with respect to Saywer-type testing condition. Even if we make it, it is generally very hard to verify Saywer-type testing condition in practice. This leads us to quest some sufficient conditions for two-weight norm inequalities of $\M_\alpha$ and $\I_\alpha$. This kind of conditions should mainly concerned with $A_p$ like conditions.

In this paper, we are mainly concerned with $A_p$ like conditions that are sufficient for two-weight norm inequalities of $\M_\alpha$ and $\I_\alpha$. First, we will work with the multiple version of entropy conditions and try to obtain the entropy bounds of $\M_\alpha$ and $\I_\alpha$. For simplicity, we only give the results and the proofs in the case m=2, although our results still hold for general $m\ge 2$.
\subsection{Main results}\label{Sec-Def-main}
First, we give one definition related to multiple weights.
\begin{definition}[\textbf{Multiple weights class}]
Let $0 \leq \alpha < mn$, $\frac{1}{p}=\frac{1}{p_{1}}+\cdots+\frac{1}{p_{m}}$ with $1 < p_{1},\cdots,p_m < \infty$,
and $0 < p \leq q < \infty$. Let $w,\sigma_i(i=1,\ldots,m)$ be nonnegative and locally integrable functions on $\Rn$,
and $\nu_{\vec{\sigma}}=\prod_{i=1}^m\sigma_i^{p/{p_i}}$. We define

\begin{eqnarray*}
[w,\vec{\sigma}]_{A_{(\vec{p},q)}}
&:=&\sup_{Q} \mathcal{A}_{\vec{p},q}(w,\vec{\sigma};Q) < \infty, \\
{}[w,\vec{\sigma}]_{A_{(\vec{p},q)}A_{\infty}^{\exp}}
&:=& \sup_{Q} \mathcal{A}_{\vec{p},q}(w,\vec{\sigma};Q) A^{\exp}_{\infty}(\nu_{\vec{\sigma}};Q)^{\frac1p} <\infty, \\
{}[w,\vec{\sigma}]_{A_{(\vec{p},q)H_{\vec{p}}^{\infty}}}
&:=& \sup_{Q} \mathcal{A}_{\vec{p},q}(w,\vec{\sigma};Q) \prod_{i=1}^m A_{\infty}^{\exp}(\sigma_i;Q)^{\frac{1}{p_i}}
<\infty,
\end{eqnarray*}
where
$$
A_{\infty}^{\exp}(w;Q)
:= \langle w \rangle_Q \exp \big(\langle \log w^{-1} \rangle_Q \big), \ \
\mathcal{A}_{\vec{p},q}(w,\vec{\sigma};Q)
:= |Q|^{\frac{1}{q}-\frac{1}{p}+\frac{\alpha}{n}} \langle w \rangle_Q^{\frac1q}
\prod_{i=1}^m \langle \sigma_{i}\rangle_Q^{\frac{1}{p_i'}}.
$$
\end{definition}

\begin{remark}
Denote
$$
[\vec{\sigma}]_{H_{\vec{p}}^\infty}
:=\sup_{Q} \prod_{i=1}^m A_{\infty}^{\exp}(\sigma_i;Q)^{\frac{p}{p_i}},\ \
[\vec{\sigma}]_{RH_{\vec{p}}} := \sup_{Q} \nu_{\vec{\sigma}}(Q)^{-1} \prod_{i=1}^m \sigma_i(Q)^{\frac{p}{p_i}},
$$
Then, it is easy to check that
\begin{align*}
&[\vec{\sigma}]_{H_{\vec{p}}^\infty} \leq [\vec{\sigma}]_{RH_{\vec{p}}} [\nu_{\vec{\sigma}}]_{A_\infty^{\exp}},\ \
[w,\vec{\sigma}]_{A_{(\vec{p},q)A_{\infty}^{\exp}}} \leq [w,\vec{\sigma}]_{A_{(\vec{p},q)}}[\nu_{\vec{\sigma}}]_{A_{\infty}^{\exp}}^{1/p}, \\
&{}[w,\vec{\sigma}]_{A_{(\vec{p},q)} A_{\infty}^{\exp}}
\leq [w,\vec{\sigma}]_{A_{(\vec{p},q)} H_{\vec{p}}^\infty}
\leq [\vec{\sigma}]_{RH_{\vec{p}}}^{1/p} [w,\vec{\sigma}]_{A_{(\vec{p},q)} A_{\infty}^{\exp}}.
\end{align*}
\end{remark}
Now, we give the definition of multilinear version of entropy conditions.
\begin{definition}[\textbf{Multilinear version of entropy conditions}]
Let $0 \leq \alpha < mn$, $\frac{1}{p}=\frac{1}{p_{1}}+\cdots+\frac{1}{p_{m}}$ with $1 < p_{1},\cdots,p_m < \infty$,
and $0 < p \leq q < \infty$. Let $w,\sigma_i(i=1,\ldots,m)$ be nonnegative and locally integrable functions on $\Rn$. We define
\begin{eqnarray*}
\lceil w,\vec{\sigma} \rceil_{\vec{p},q,\epsilon}
&:=& \sup_{Q} \mathcal{A}_{\vec{p},q}(w,\vec{\sigma};Q) \rho_{\nu_{\vec{\sigma}}}(Q)^{\frac{1}{p}} \epsilon(\rho_{\nu_{\vec{\sigma}}}(Q)) \\
\lfloor w,\vec{\sigma} \rfloor_{\vec{p},q,\vec{\epsilon},\eta}
&:=& \sup_{Q} \mathcal{A}_{\vec{p},q}(w,\vec{\sigma};Q) \rho_{w,\eta}(Q)^{\frac{1}{q'}}
\prod_{i=1}^m \rho_{\sigma_i,\epsilon_i}(Q)^{\frac{1}{p_i}} \\
{}[[ \vec{\sigma}]]_{(i,j,k),\epsilon_i}
&:=& \sup_{Q} \Big(|Q|^{\frac{\alpha}{n}} \prod_{i=1}^2 \langle \sigma_i \rangle_Q\Big)^{\frac{p_k'}{p_{ij}'}} \cdot \langle \sigma_3 \rangle_Q \gamma_{(i,j,k)}(Q) \epsilon_i(\gamma_{(i,j,k)}(Q)).
\end{eqnarray*}
where $\epsilon,\eta,\epsilon_i$ are monotonic increasing functions on $(1,\infty)$, and
\begin{eqnarray*}
\rho_{w}(Q)&:=&\frac{\int_Q M(\mathbf{1}_Q w)(x)dx}{w(Q)}  \ \text{and} \
\rho_{w,\epsilon}(Q)£º= \rho_{w}(Q) \epsilon(\rho_{w}(Q)), \\
\gamma_{(i,j,k)}(Q)
&:=& \frac{\int_Q \mathcal{M}_\alpha(\mathbf{1}_Q \sigma_i, \mathbf{1}_Q \sigma_j)(x)^{\frac{p_k'}{p_{ij}}}dx}
{\big(\int_{Q}\sigma_i^{\frac{p_{ij}}{p_i}} \sigma_j^{\frac{p_{ij}}{p_j}} dx\big)^{\frac{p_k'}{p_{ij}}}} \ \text{and} \
\frac{1}{p_{ij}}=\frac{1}{p_i}+\frac{1}{p_j}.
\end{eqnarray*}
\end{definition}

The main results of this paper are as follows:

\begin{theorem}\label{Theorem M-alpha}
Let $0 \leq \alpha < 2n$, $0 < p \leq q < \infty$ and $\frac{1}{p}=\frac{1}{p_1}+\frac{1}{p_2}$ with $1 < p_1,p_2 < \infty$. Suppose that $\sigma_1,\sigma_2,w$ are weights on $\Rn$. Let $\epsilon$ be a monotonic increasing function on $(1,\infty)$ that satisfies $\int_{1}^\infty \frac{dt}{t \epsilon(t)^q} < \infty$. Then, the following inequality holds
$$
\big\| \mathcal{M}_{\alpha}(f_1 \sigma_1, f_2 \sigma_2) \big\|_{L^q(w)}
\lesssim \lceil w,\vec{\sigma} \rceil_{\vec{p},q,\epsilon} \prod_{i=1}^2 \big\| f_i \big\|_{L^{p_i}(\sigma_i)}.
$$
\end{theorem}

\begin{theorem}\label{Theorem I-alpha-1}
Let $0 \leq \alpha < 2n$, $0 < p \leq q < \infty$ and $\frac{1}{p}=\frac{1}{p_1}+\frac{1}{p_2}$ with $1 < p_1,p_2,q < \infty$. Suppose that $\sigma_1,\sigma_2,w$ are weights on $\Rn$. Let $\epsilon_1,\epsilon_2,\eta$ be monotonic increasing functions on $(1,\infty)$ that satisfies $\int_{1}^\infty \frac{dt}{t \epsilon_i(t)^{p_i}} < \infty$ and $\int_{1}^\infty \frac{dt}{t \eta(t)^{q'}} < \infty$. Then, the following inequality holds
$$
\big\| \mathcal{I}_{\alpha}(f_1 \sigma_1, f_2 \sigma_2) \big\|_{L^q(w)}
\lesssim \lfloor w,\vec{\sigma} \rfloor_{\vec{p},q,\vec{\epsilon},\eta}
\prod_{i=1}^2 \big\| f_i \big\|_{L^{p_i}(\sigma_i)}.
$$
\end{theorem}

\begin{theorem}\label{Theorem I-alpha-2}
Let $0 \leq \alpha < 2n$, and $1 < p_i < \infty\ (i=1,2,3)$ satisfying $\frac{1}{p_i} + \frac{1}{p_j} \geq 1$ for $i \neq j$.
Suppose that $\sigma_1,\sigma_2,\sigma_3$ are weights on $\Rn$. Let $\epsilon_i$ be monotonic increasing functions on $(1,\infty)$ such that $\int_{1}^\infty \frac{dt}{t \epsilon_i(t)^{1/{p_i'}}} < \infty$, $i=1,2,3$. Then, the following inequality holds
\begin{eqnarray}
\big\| \I_{\alpha}(f_1 \sigma_1, f_2 \sigma_2) \big\|_{L^{p_3'}(\sigma_3)}
& \lesssim & \sum_{(i,j,k) \in \Omega} [[ \vec{\sigma}]]_{(i,j,k),\epsilon_i}^{1/{p_k'}}
\prod_{i=1}^2 \big\| f_i \big\|_{L^{p_i}(\sigma_i)}; \\
\big\| \I_{\alpha}(f_1 \sigma_1, f_2 \sigma_2) \big\|_{L^{p_3',\infty}(\sigma_3)}
& \lesssim & \sum_{\substack{i \neq 3 \\ (i,j,k) \in \Omega}} [[ \vec{\sigma}]]_{(i,j,k),\epsilon_i}^{1/{p_k'}}
\prod_{i=1}^2 \big\| f_i \big\|_{L^{p_i}(\sigma_i)};
\end{eqnarray}
where $\Omega$ is the set of all permutations of $(1, 2, 3)$.
\end{theorem}

The article is organized as follows: In Section $\ref{Sec-2}$, some notations and lemmas will be given. In Section $\ref{Sec-3}$, we will demonstrate Theorem $\ref{Theorem M-alpha}$ and Theorem $\ref{Theorem I-alpha-1}$. Section $\ref{Sec-4}$ will be devoted to complete the proofs of Theorem $\ref{Theorem I-alpha-2}$.

\section{Preliminaries}\label{Sec-2}
First, we present some definitions and lemmas, which will be used later.
\begin{definition}
A collection, $\D$ of cubes is said to be a dyadic grid if it satisfies
\begin{enumerate}
\item [(1)] The side length of every $Q \in \D$ equals $2^k$ for some $k \in \Z$.
\item [(2)] For any $Q,R \in \D$, $Q \cap R = \{Q,R,\emptyset\}$.
\item [(3)] $\Rn = \bigcup_{Q \in \D_k} Q$, $\D_k=\{Q \in \D; \ell(Q)=2^k\}$ for any $k \in \Z$.
\end{enumerate}
\end{definition}

\begin{definition}
A subset $\S$ of a dyadic grid is said to be spare, if for every $Q \in \S$ there holds that
$$
\Big| \bigcup_{\substack{Q' \in \S \\ Q' \subsetneq Q}} Q' \Big| \leq \frac12 |Q|.
$$
Equivalently, if $E(Q) = Q \setminus \bigcup_{\substack{Q' \in \S \\ Q' \subsetneq Q}} Q'$, then
the sets $\{E(Q)\}_{Q \in \S}$ are pairwise disjoint and $|Q| \leq 2 |E(Q)|$.
\end{definition}
\begin{definition}
Let $0 \leq \alpha < mn$ and $\D,\S$ be a given dyadic grid and a spare set. The dyadic versions of multilinear fractional maximal and fractional integral operators are defined by
\begin{eqnarray*}
\M_{\alpha}^{\D}(\vec{f})(x)
&:=& \sup_{Q \in \D} |Q|^{\frac{\alpha}{n}} \prod_{i=1}^m \langle |f_i| \rangle_Q \cdot \mathbf{1}_Q(x),\\
\I_{\alpha}^{\D}(\vec{f})(x)
&:=& \sum_{Q \in \D} |Q|^{\frac{\alpha}{n}} \prod_{i=1}^m \langle f_i \rangle_Q \cdot \mathbf{1}_Q(x), \\
\T_{\S,\alpha}(\vec{f})(x)
&:=& \sum_{Q \in \S} |Q|^{\frac{\alpha}{n}}\prod_{i=1}^m \langle f_i \rangle_Q \cdot \mathbf{1}_{E(Q)}(x).
\end{eqnarray*}
\end{definition}
We will need the following lemma given by Hyt\"{o}nen and P\'{e}rez in \cite{HP}.
\begin{lemma}
There are $2^n$ dyadic grids $\D_t$, $t \in \{0, 1/3\}^n$ such that for any cube $Q \subset \Rn$ there exists a cube $Q_t \in \D_t$ satisfying $Q \subset Q_t$ and $\ell(Q_t) \leq 6 \ell(Q)$, where the dyadic grid $\D_t$ is defined by
$$
\D_t := \big\{2^{-k}([0, 1)^n + m + (-1)^kt) : k \in \Z,m \in \Z^n \big\},\ \ t \in \{0, 1/3\}^n.
$$
\end{lemma}

We also need the following lemma.
\begin{lemma}\label{pointwise}
Let $\D$ be a dyadic grid. For any non-negative integrable $f_i (i=1,\ldots, m)$, there exist sparse families $\S \subset \D$ such that for all $x \in \Rn$, it holds that
\begin{eqnarray}
\label{M-alpha} \M_{\alpha}(\vec{f})(x) \simeq \sup_{t \in \{0,1/3\}^n}\M_{\alpha}^{\D_t}(\vec{f})(x),&{}&
\M_{\alpha}^{\D}(\vec{f})(x) \simeq \T_{\S,\alpha}(\vec{f})(x); \\
\label{I-alpha} \I_{\alpha}(\vec{f})(x) \simeq \sup_{t \in \{0,1/3\}^n}\I_{\alpha}^{\D_t}(\vec{f})(x),&{}&
\I_{\alpha}^{\D}(\vec{f})(x) \simeq \T_{\S,\alpha}(\vec{f})(x).
\end{eqnarray}
\end{lemma}
The proof of $(\ref{M-alpha})$ can be found in \cite{LMS} and $(\ref{I-alpha})$ was shown in \cite{LS-2}.

We will need to apply the following multilinear version of Carleson embedding theorem \cite{S} at certain key points in the proofs of our results.
\begin{lemma}[\textbf{Carleson embedding theorem}]\label{Carleson Embedding}
Let $0 < p \leq q < \infty$ and $\frac{1}{p} = \frac{1}{p_1} + \cdots + \frac{1}{p_m}$ satisfying $1 < p_1,\ldots,p_m < \infty$.
Suppose that the nonnegative numbers $\{ c_Q \}_{Q}$ satisfy
$$
\sum_{Q \subset Q'} c_Q \leq A \ \nu_{\vec{\sigma}}(Q')^{q/p} , \text{ for any } Q' \in \D,
$$
where $\sigma_i \ (i = 1, \cdots, m)$ are weights and $\nu_{\vec{\sigma}} = \prod_{i=1}^m \sigma_i^{p/{p_i}}$.
Then for all nonnegative functions $f_i \in L^{p_i}(\sigma_i)$, we have
\begin{align*}
\sum_{Q \in \D} c_Q \prod_{i=1}^m \big( \langle f_i \rangle_Q^{\sigma_i} \big)^q
\lesssim A \big\| \mathcal{M}_{\vec{\sigma}}^d(\vec{f}) \big\|_{L^{p,q}(\nu_{\vec{\sigma}})}^q
\lesssim A \prod_{i=1}^m \big\| f_i \big\|_{L^{p_i}(\sigma_i)}^q,
\end{align*}
where $L^{p,q}(w)$ is the Lorentz space defined by
$$ \big\| f \big\|_{L^{p,q}(w)} = \bigg[ \int_{0}^\infty \Big( \lambda w \big(\{ x \in \Rn; |f(x)| > \lambda \}\big)^{1/p} \Big)^q \frac{d\lambda}{\lambda} \bigg]^{1/q} < \infty.$$
\end{lemma}


\section{Proofs of Theorems $\ref{Theorem M-alpha}$-$\ref{Theorem I-alpha-1}$}\label{Sec-3}
In this section, our aim is to demonstrate Theorem $\ref{Theorem M-alpha}$ and Theorem $\ref{Theorem I-alpha-1}$
by making use of dyadic techniques (see for examples, \cite{H-A2} and \cite{Lacey-S}).
\subsection{Proof of Theorem $\ref{Theorem M-alpha}$.}
Let $\S$ be any sparse set of $\D$. By Lemma $\ref{pointwise}$, it suffices to show that
\begin{equation}\label{L-S-alpha}
\big\| \T_{\S,\alpha}(f_1 \sigma_1, f_2 \sigma_2) \big\|_{L^q(w)}
\lesssim  \lceil w,\vec{\sigma} \rceil_{\vec{p},q,\epsilon}
\prod_{i=1}^2 \big\| f_i \big\|_{L^{p_i}(\sigma_i)}.
\end{equation}
We may assume that each $f_i$ is a non-negative function for $i=1,2$.
Denote
$$
\S_k := \Big\{Q \in \S;\ 2^{-k} \lceil w,\vec{\sigma}\rceil_{\vec{p},q,\epsilon} \leq \Gamma(Q) \leq 2^{-k+1} \lceil w,\vec{\sigma}\rceil_{\vec{p},q,\epsilon} \Big\},
$$
where
$
\Gamma(Q):= |Q|^{\frac{1}{q}-\frac{1}{p}+\frac{\alpha}{n}} \langle w \rangle_Q^{\frac1q}
\prod_{i=1}^2 \langle \sigma_{i}\rangle_Q^{\frac{1}{p_i'}} \cdot \rho_{\nu_{\vec{\sigma}}}(Q)^{\frac{1}{p}}\epsilon(\rho_{\nu_{\vec{\sigma}}}(Q)).
$
Using the pairwise disjointness of the sets $\{E(Q)\}_{Q \in \S}$, we deduce that
\begin{align*}
\big\| \T_{\S,\alpha}(f_1 \sigma_1, f_2 \sigma_2) \big\|_{L^q(w)}^q
&= \sum_{k=1}^\infty \sum_{Q \in \S_k} \Big(|Q|^{\frac{\alpha}{n}} \prod_{i=1}^2 \langle f_i \sigma_i \rangle_Q \Big)^q w(E(Q))
:= \sum_{k=1}^\infty \Delta_k.
\end{align*}
To obtain the bound of $\Delta_k$, we need to introduce the notion $$
c_Q = \big(|Q|^{\alpha/n} \langle \sigma_1 \rangle_Q \langle \sigma_2 \rangle_Q \big)^q w(E(Q)).
$$Then, it is easy to see that
\begin{align*}
\Delta_k
=\sum_{Q \in \S_k} c_Q \Big(\langle f_1 \rangle_Q^{\sigma_1} \langle f_2 \rangle_Q^{\sigma_2} \Big)^q.
\end{align*}
In order to apply the Carleson embedding theorem, we need to analyze $\{c_Q\}_{Q \in \S_k}$.
Fix $Q' \in \S_k$. Since $\Gamma(Q) \simeq 2^{-k}\lceil w,\vec{\sigma} \rceil_{\vec{p},q,\epsilon}$ for each $Q \in \S_k$, we get
\begin{align*}
\sum_{Q \in \S_k: Q \subset Q'}c_Q
&\leq \sum_{Q \in \S_k: Q \subset Q'} \Gamma(Q)^q
\bigg( \frac{\sigma_1(Q)^{\frac{1}{p_1}} \sigma_2(Q)^{\frac{1}{p_2}}}
{\rho_{\nu_{\vec{\sigma}}}(Q)^{\frac{1}{p}}\epsilon(\rho_{\nu_{\vec{\sigma}}}(Q))} \bigg)^q \\
&\lesssim 2^{-kq}\lceil w,\vec{\sigma} \rceil_{\vec{p},q,\epsilon}^q \sum_{Q \in \S_k: Q \subset Q'}
\frac{\nu_{\vec{\sigma}}(Q)^{q/p}}{\rho_{\nu_{\vec{\sigma}}}(Q)^{q/p}\epsilon(\rho_{\nu_{\vec{\sigma}}}(Q))^q} \\
&:= 2^{-kq}\lceil w,\vec{\sigma} \rceil_{\vec{p},q,\epsilon}^q \Delta_k'.
\end{align*}
Now, we tentatively claim that
\begin{equation}\label{Delta k'}
\Delta_k' \lesssim \nu_{\vec{\sigma}}(Q')^{q/p}.
\end{equation}
Therefore, if the above claim is true, we actually obtain that
$$
\sum_{Q \in \S_k: Q \subset Q'}c_Q
\lesssim 2^{-kq}\lceil w,\vec{\sigma} \rceil_{\vec{p},q,\epsilon}^q \nu_{\vec{\sigma}}(Q')^{q/p},
$$
and
$$
\Delta_k \lesssim
2^{-kq}\lceil w,\vec{\sigma} \rceil_{\vec{p},q,\epsilon}^q \prod_{i=1}^2 \big\| f_i \big\|_{L^{p_i}(\sigma_i)}^q
$$
provided by Lemma $\ref{Carleson Embedding}$. Therefore, it yields that
\begin{equation*}
\big\| \T_{\S,\alpha}(\vec{f} \cdot \vec{\sigma}) \big\|_{L^q(w)}
=\sum_{k=1}^\infty \Delta_k
\lesssim  \lceil w,\vec{\sigma} \rceil
\prod_{i=1}^2 \big\| f_i \big\|_{L^{p_i}(\sigma_i)}.
\end{equation*}
This shows that inequality $(\ref{L-S-alpha})$ is true.

Now, we are in the position to prove $(\ref{Delta k'})$.
Set
$$
\S_{k,j}:=\big\{Q \in \S_k;\ Q \subset Q',\ 2^{j-1} \leq \rho_{\nu_{\vec{\sigma}}}(Q) < 2^j \big\},
$$
and $\S_{k,j}^*$ is the collection of maximal elements in $\S_{k,j}$. Thereby, we have
\begin{align*}
\Big(\sum_{Q \in \S_{k,j}} \nu_{\vec{\sigma}}(Q)^{q/p}\Big)^{p/q}
&\leq \sum_{Q^* \in \S_{k,j}^*} \sum_{Q \subset Q^*} \nu_{\vec{\sigma}}(Q) \\
&\lesssim \sum_{Q^* \in \S_{k,j}^*} \sum_{Q \subset Q^*} \int_{E(Q)} \langle \mathbf{1}_{Q^*}\nu_{\vec{\sigma}} \rangle_Q \mathbf{1}_Q(x)dx \\
&\lesssim \sum_{Q^* \in \S_{k,j}^*} \sum_{Q \subset Q^*} \int_{E(Q)} \sup_{P \in \D}\langle \mathbf{1}_{Q^*}\nu_{\vec{\sigma}} \rangle_P \mathbf{1}_P(x)dx \\
&\leq \sum_{Q^* \in \S_{k,j}^*}\int_{Q^*} \sup_{P \in \D}\langle \mathbf{1}_{Q^*}\nu_{\vec{\sigma}} \rangle_P \mathbf{1}_P(x)dx \\
&\leq \sum_{Q^* \in \S_{k,j}^*}\int_{Q^*} M(\mathbf{1}_{Q^*}\nu_{\vec{\sigma}})(x)dx \\
&=\sum_{Q^* \in \S_{k,j}^*}\nu_{\vec{\sigma}}(Q^*) \rho_{\nu_{\vec{\sigma}}}(Q^*)
\lesssim 2^j \nu_{\vec{\sigma}}(Q').
\end{align*}
Consequently, we deduce that
\begin{align*}
\Delta_k'
&\leq \sum_{j=0}^\infty \frac{1}{2^{jq/p} \epsilon(2^j)^q}\sum_{Q \in \S_{k,j}} \nu_{\vec{\sigma}}(Q)^{q/p} \lesssim \nu_{\vec{\sigma}}(Q')^{q/p} \sum_{j=0}^\infty \frac{1}{\epsilon(2^j)^q}\\&\lesssim \nu_{\vec{\sigma}}(Q')^{q/p} \int_{1}^\infty \frac{dt}{t \epsilon(t)^q}
\lesssim \nu_{\vec{\sigma}}(Q')^{q/p}.
\end{align*}
The proof of $(\ref{Delta k'})$ is finished.

\qed
\subsection{Proof of Theorem $\ref{Theorem I-alpha-1}$.} By duality, we have
\begin{align*}
\big\| \T_{\S,\alpha}(\vec{f}\cdot\vec{\sigma}) \big\|_{L^q(w)}
&=\sup_{||g||_{L^{q'}(w) \leq 1}} \bigg|\sum_{Q \in \S}\big( |Q|^{\frac{\alpha}{n}}\langle f_1 \sigma_1 \rangle_Q \langle f_2 \sigma_2 \rangle_Q \big) \int_Q g(x) w dx \bigg| \\
&:= \sup_{||g||_{L^{q'}(w) \leq 1}} \big| \mathfrak{I}(g) \big|.
\end{align*}
Denote
\begin{eqnarray*}
\S_k = \big\{Q \in \S;\ 2^k < \lambda_Q \leq 2^{k+1} \big\},\ \
\lambda_Q = \mathcal{A}_{\vec{p},q}(w,\vec{\sigma};Q) \rho_{w,\eta}(Q)^{\frac{1}{q'}}
\prod_{i=1}^2 \rho_{\sigma_i,\epsilon_i}(Q)^{\frac{1}{p_i}}.
\end{eqnarray*}
Then, we have $k \leq K_0 := \log_2 \lfloor w,\vec{\sigma} \rfloor_{\vec{p},q,\vec{\epsilon},\eta}$. Therefore, by the H\"{o}lder inequality, it now follows that
\begin{align*}
\mathfrak{I}(g)
&=\sum_{k=1}^{K_0} \sum_{Q \in \S_k} \lambda_Q \prod_{i=1}^2 \frac{\langle f_i \rangle_Q^{\sigma_i} \sigma_i(Q)^{\frac{1}{p_i}}}{\rho_{\sigma_i,\epsilon_i}(Q)^{\frac{1}{p_i}}}
\frac{\langle g \rangle_Q^w w(Q)^{\frac{1}{q'}}}{\rho_{w,\eta}(Q)^{\frac{1}{q'}}} \\
&\lesssim \sum_{k=1}^{K_0} 2^{k} \bigg(\sum_{Q \in \S_k} \prod_{i=1}^2 \frac{ \big(\langle f_i \rangle_Q^{\sigma_i}\big)^q \sigma_i(Q)^{\frac{q}{p_i}}}{\rho_{\sigma_i,\epsilon_i}(Q)^{\frac{q}{p_i}}}\bigg)^{\frac{1}{q}}
\bigg(\sum_{Q \in \S_k} \big(\langle g \rangle_Q^w\big)^{q'} \frac{w(Q)}{\rho_{w,\eta}(Q)}\bigg)^{\frac1{q'}} \\
&\lesssim \lfloor w,\vec{\sigma} \rfloor_{\vec{p},q,\vec{\epsilon},\eta} \prod_{i=1}^2 \bigg(\sum_{Q \in \S_k}  \big(\langle f_i \rangle_Q^{\sigma_i}\big)^{p_i} \frac{\sigma_i(Q)}{\rho_{\sigma_i,\epsilon_i}(Q)} \bigg)^{\frac{1}{p_i}}
\bigg(\sum_{Q \in \S_k} \big(\langle g \rangle_Q^w\big)^{q'} \frac{w(Q)}{\rho_{w,\eta}(Q)}\bigg)^{\frac1{q'}}.
\end{align*}
By the Carleson embedding theorem $\ref{Carleson Embedding}$, it is enough to show that for each $Q' \in \S_k$
$$
\sum_{\substack{Q \in \S_k \\ Q \subset Q'}}\frac{\sigma_i(Q)}{\rho_{\sigma_i,\epsilon_i}(Q)} \lesssim \sigma_i(Q'),\ \ i=1,2,3,
$$
where $\sigma_3=w$ and $\epsilon_3=\eta$. A completely analogous calculation to that of the preceding subsection yields the desired results.

\qed

\section{Proof of Theorem $\ref{Theorem I-alpha-2}$} \label{Sec-4}
In this section, we shall give the proof of Theorem $\ref{Theorem I-alpha-2}$. We need the following two-weight characterization of $\I_{\alpha}^{\S}$, which was proved in \cite{HHL} and \cite{LS-2}.
\begin{lemma}\label{Characterization}
Let $\D$ be a dyadic grid and $\S \subset \D$ be a sparse family. Suppose that $\sigma_1, \sigma_2,\sigma_3$ are positive Borel measures and $1 < p_i < \infty\ (i=1,2,3)$ satisfying $\frac{1}{p_i} + \frac{1}{p_j} \geq 1$ for $i \neq j$. Then
\begin{enumerate}
\item [(1)] The strong type inequality 
$$
\big\| \mathcal{I}_{\alpha}^{\mathcal{S}}(f_1 \sigma_1, f_2 \sigma_2)\big\|_{L^q(\sigma_3)}
\leq \mathfrak{N}\prod_{i=1}^2 \big\| f_i \big\|_{L^{p_i}(\sigma_i)}
$$
holds if and only if the following test conditions hold for any triple $(i,j,k) \in \Omega$,
\begin{eqnarray*}
\mathfrak{T}_{\S,(i,j,k)} &:=& \sup_{R \in \S} \frac{
\Big\| \sum_{\substack{Q \in \S \\ Q \subset R}}|Q|^{\frac{\alpha}{n}}\langle \sigma_j \rangle_Q \langle \sigma_k \rangle_Q \mathbf{1}_Q \Big\|_{L^{p_i'}(\sigma_i)}}{\sigma_j(R)^{1/{p_j}} \sigma_k(R)^{1/{p_k}}} < \infty.
\end{eqnarray*}
\item [(2)] The weak type inequality 
$$
\big\| \mathcal{I}_{\alpha}^{\mathcal{S}}(f_1 \sigma_1, f_2 \sigma_2)\big\|_{L^{q,\infty}(\sigma_3)}
\leq \mathfrak{N}_{weak}\prod_{i=1}^2 \big\| f_i \big\|_{L^{p_i}(\sigma_i)}
$$
holds if and only if $\mathfrak{T}_{\S,(i,j,k)} < \infty$,
for any triple $(i,j,k) \in \Omega$ and $i \neq 3$.
\end{enumerate}
Moreover, the best constants satisfy
$\mathfrak{N} \simeq \sum_{(i,j,k) \in \Omega}\mathfrak{T}_{\S,(i,j,k)}$,
$\mathfrak{N}_{weak} \simeq \sum_{i \neq 3,(i,j,k) \in \Omega}\mathfrak{T}_{\S,(i,j,k)}$.
\end{lemma}
\noindent\textbf{Proof of Theorem $\ref{Theorem I-alpha-2}$.}
By Lemma $\ref{Characterization}$, it suffices to show
$$
\mathfrak{T}_{\S,(i,j,k)} \lesssim [[ \vec{\sigma}]]_{(i,j,k),\epsilon_i},\ \text{ for each } (i,j,k) \in \Omega.
$$
By symmetry, we only focus on estimating the case $(i,j,k)=(1,2,3)$. For convenience, we write $q=p_3'$, $p=p_{12}$ and $\gamma=\gamma_{(1,2,3)}$. From now on, we fix the cube $R \in \D$ and introduce the notations
\begin{eqnarray*}
\mathscr{A}(R) &:=& \Big\| \sum_{Q \in \S}|Q|^{\frac{\alpha}{n}}\langle \sigma_1 \rangle_Q \langle \sigma_2 \rangle_Q
\mathbf{1}_Q \Big\|_{L^{q}(R,w)},\\
\mathscr{B}(Q) &:=& \Big(|Q|^{\frac{\alpha}{n}} \prod_{i=1}^m \langle \sigma_i \rangle_Q\Big)^{\frac{q}{p'}} \cdot
\langle w \rangle_Q \gamma(Q) \epsilon_1(\gamma(Q)).
\end{eqnarray*}
Then, we make a partition of $\S$ by setting
$$
\S_{a,b}:=\big\{Q \in \S;\ Q \subset R, 2^a < \mathscr{B}(Q) \leq 2^{a+1}, 2^b < \gamma(Q) \leq 2^{b+1}\big\}.
$$
Note that $2^a \leq [[ \vec{\sigma}]]_{(1,2,3),\epsilon_1}$. Now we construct the stopping cubes $\F$. Let $\F$
be the minimal subset of $\S_{a,b}$ containing the maximal cubes in $\S_{a,b}$ such that whenever $F \in \F$, the maximal cubes
$Q \subset F$, $Q \in \S_{a,b}$ with
$|Q|^{\frac{\alpha}{n}} \langle \sigma_1 \rangle_Q \langle \sigma_1 \rangle_Q > 4|F|^{\frac{\alpha}{n}} \langle \sigma_1 \rangle_F \langle \sigma_1 \rangle_F$ are also in $\F$. Denote by $\pi_{\F}(Q)$ the minimal cube in $\F$ which contains $Q$. Denote
$$
\S_{a,b}^k := \Big\{Q \in \S_{a,b};\ |Q|^{\frac{\alpha}{n}} \langle \sigma_1 \rangle_Q \langle \sigma_2 \rangle_Q
\simeq  2^{-k} |\pi_{\F}(Q)|^{\frac{\alpha}{n}} \langle \sigma_1 \rangle_{\pi_{\F}(Q)} \langle \sigma_2 \rangle_{\pi_{\F}(Q)} \Big\}.
$$
Then the Minkowski inequality implies that
\begin{equation}\label{A}
\mathscr{A}(R)
\leq \sum_{a,b} \sum_{k=1}^\infty \Big\| \sum_{Q \in \S_{a,b}^k}|Q|^{\frac{\alpha}{n}}\langle \sigma_1 \rangle_Q \langle \sigma_2 \rangle_Q \mathbf{1}_Q \Big\|_{L^q(w)}
:= \sum_{a,b} \sum_{k=1}^\infty \Theta_{a,b}^k.
\end{equation}
For each $F \in \F$, write
$$
\Psi_F := \sum_{\substack{Q \in \S_{a,b}^k \\ \pi_{\F}(Q)=F}} |Q|^{\frac{\alpha}{n}}\langle \sigma_1 \rangle_Q \langle \sigma_2 \rangle_Q,\
\Psi_{F,j}:=\Psi_F \mathbf{1}_{\{\Psi_F \simeq j 2^{-k} |F|^{\frac{\alpha}{n}}\langle \sigma_1 \rangle_F \langle \sigma_2 \rangle_F \}}.
$$
Applying the H\"{o}lder inequality, we may obtain that
\begin{equation}\aligned\label{a-b-k}
\Theta_{a,b}^k
&\leq \Big\| \Big(\sum_{j=1}^\infty j^{-\frac{2}{q'}} j^{\frac{2}{q'}}\sum_{F \in \mathcal{F}} \Psi_{F,j} \Big)^q \Big\|_{L^1(w)}^{\frac{1}{q}} \\
&\leq \Big(\sum_{j=1}^\infty j^{-\frac{2}{q'}q'}\Big)^{\frac1q} \Big\| \sum_{j=1}^\infty j^{\frac{2q}{q'}}
\sum_{F \in \F} \Psi_{F,j}^q \Big\|_{L^1(w)}^{\frac{1}{q}} \\
&\lesssim \bigg(\sum_{j=1}^\infty j^{2(q-1)} \sum_{F \in \F} \int_{Q_0} \Psi_{F,j}(x)^q w \ dx \bigg)^{\frac{1}{q}}.
\endaligned
\end{equation}
Therefore, we are in a position to consider the contribution of the integral in the above inequality. Before doing that,
we first claim that the following estimate is true:
\begin{equation}\label{Psi}
w\big(\{x;\Psi_F(x) > \lambda 2^{-k} |F|^{\frac{\alpha}{n}}\langle \sigma_1 \rangle_F \langle \sigma_2 \rangle_F \}\big)
\lesssim 2^{-\lambda} w(F).
\end{equation}
By $(\ref{Psi})$ and noticing the fact that the set
$\{x;\Psi_F(x) > \lambda 2^{-k} |F|^{\frac{\alpha}{n}}\langle \sigma_1 \rangle_F \langle \sigma_2 \rangle_F \}$
coincides with $F$ if $0 < \lambda < j/2$ and is empty if $\lambda > j$, it is easy to obtain that
$$
\int_{0}^\infty q \lambda^{q-1}w\big(\{x;\Psi_F(x) > \lambda 2^{-k} |F|^{\frac{\alpha}{n}}\langle \sigma_1 \rangle_F \langle \sigma_2 \rangle_F \}\big) d\lambda
\lesssim j^q 2^{-\frac{j}{2}}w(F).
$$
Hence, it now follows that
\begin{equation}\label{Psi-F-j}
\int_{Q_0} \Psi_{F,j}(x)^q w \ dx
\lesssim 2^{-kq}|F|^{\frac{\alpha q}{n}}\langle \sigma_1 \rangle_F^q \langle \sigma_2 \rangle_F^q \big(j^q 2^{-\frac{j}{2}}w(F)\big).
\end{equation}
By estimates $(\ref{a-b-k})$ and $(\ref{Psi-F-j})$, we get
\begin{align*}
\Big(\sum_{k=1}^\infty \Theta_{a,b}^k \Big)^q
&\lesssim \Big(\sum_{k=1}^\infty 2^{-k} \Big)^q \sum_{j=1}^\infty j^{3q-2} 2^{-\frac{j}{2}} \sum_{F \in \F}
|F|^{\frac{\alpha q}{n}}\langle \sigma_1 \rangle_F^q \langle \sigma_2 \rangle_F^q w(F) \\
&\lesssim \frac{2^a}{2^b \epsilon_1(2^b)}\sum_{F \in \F} \Big(|F|^{\frac{\alpha}{n}}\langle \sigma_1 \rangle_F \langle \sigma_2 \rangle_F\Big)^{\frac{q}{p}} |F|.
\end{align*}
Let $\F^*$ be the maximal elements of $\F$, then we obtain
\begin{align*}
\Big(\sum_{k=1}^\infty \Theta_{a,b}^k \Big)^q
&\lesssim \frac{2^a}{2^b \epsilon_1(2^b)}\sum_{F^* \in \F^*} \sum_{F^* \supset F \in \F} \int_{E(F)}
\mathcal{M}_{\alpha}(\mathbf{1}_{F^*} \sigma_1, \mathbf{1}_{F^*} \sigma_2)(x)^{\frac{q}{p}} \ dx \\
&\leq \frac{2^a}{2^b \epsilon_1(2^b)}\sum_{F^* \in \F^*} \int_{F^*}
\mathcal{M}_{\alpha}(\mathbf{1}_{F^*} \sigma_1, \mathbf{1}_{F^*} \sigma_2)(x)^{\frac{q}{p}} \ dx \\
&\leq \frac{2^a}{\epsilon_1(2^b)}\sum_{F^* \in \F^*} \nu_{\vec{\sigma}}(F^*)^{\frac{q}{p}}
\leq \frac{2^a}{\epsilon_1(2^b)} \Big(\sum_{F^* \in \F^*} \nu_{\vec{\sigma}}(F^*)\Big)^{\frac{q}{p}} \\
&\leq \frac{2^a}{\epsilon_1(2^b)} \nu_{\vec{\sigma}}(R)^{\frac{q}{p}}
\leq \frac{2^a}{\epsilon_1(2^b)} \big[\sigma_1(R)^{\frac{1}{p_1}} \sigma_2(R)^{\frac{1}{p_2}} \big]^q.
\end{align*}
Consequently, the inequality $(\ref{A})$ gives that
$$
\mathscr{A}(R)
\lesssim \sum_{a,b}\frac{2^{a/q}}{\epsilon_1(2^b)^{1/q}} \sigma_1(R)^{\frac{1}{p_1}} \sigma_2(R)^{\frac{1}{p_2}}
\lesssim [[ \vec{\sigma}]]_{(1,2,3),\epsilon_1}^{1/q} \int_{1}^\infty \frac{dt}{t \epsilon_1(t)^{1/q}} \sigma_1(R)^{\frac{1}{p_1}} \sigma_2(R)^{\frac{1}{p_2}}.
$$
This shows that
$$
\mathfrak{T}_{\S,(1,2,3)} \lesssim [[ \vec{\sigma}]]_{(1,2,3),\epsilon_1}^{1/{q}}.
$$

We are left to prove the claim, inequality $(\ref{Psi})$. If $w$ is the Lebesgue measure, the inequality is obvious. For any $Q \in \S_{a,b}^k$ satisfying $\pi_{\F}(Q) = F$, it holds that
$$
2^a \simeq \mathscr{B}(Q) \simeq \big(2^{-k}|F|^{\frac{\alpha}{n}}\langle \sigma_1 \rangle_F \langle \sigma_2 \rangle_F\big)^{\frac{q}{p'}} \langle w \rangle_Q 2^b \epsilon_1(2^b).
$$
Let $\S_{a,b}^{k,*}$ be the maximal cubes in $\S_{a,b}^k$ and $
\Lambda_F := \frac{2^a}{2^b \epsilon_1(2^b)} \big(2^{-k} |F|^{\frac{\alpha}{n}}\langle \sigma_1 \rangle_F \langle \sigma_2 \rangle_F\big)^{-\frac{q}{p'}},
$. Note that the set
$\{x;\Psi_F(x) > \lambda 2^{-k} |F|^{\frac{\alpha}{n}}\langle \sigma_1 \rangle_F \langle \sigma_2 \rangle_F \}$
is the union of maximal cubes $P \in \S_{a,b}^k$ with $\pi_{\F}(P)=F$ and
$\inf_{x \in P} \Psi_F(x) > \lambda 2^{-k} |F|^{\frac{\alpha}{n}}\langle \sigma_1 \rangle_F \langle \sigma_2 \rangle_F$.
Then, it yields that
\begin{align*}
&w\big(\{x;\Psi_F(x) > \lambda 2^{-k} |F|^{\frac{\alpha}{n}}\langle \sigma_1 \rangle_F \langle \sigma_2 \rangle_F \}\big)\\ 
&\simeq \Lambda_F \big|\{x;\Psi_F(x) > \lambda 2^{-k} |F|^{\frac{\alpha}{n}}\langle \sigma_1 \rangle_F \langle \sigma_2 \rangle_F \}\big| \\
&\lesssim \Lambda_F \ 2^{-\lambda} \sum_{Q^* \in \S_{a,b}^{k,*}} |Q^*|
\lesssim 2^{-\lambda} \sum_{Q^* \in \S_{a,b}^{k,*}} w(Q^*) \\
&\lesssim 2^{-\lambda} w(F).
\end{align*}
\qed

\end{document}